\title{A Tight Erd\H{o}s-Stone Bound for All Graph Densities}

\author{Asaf Shapira
\thanks{School of Mathematics, Tel Aviv University, Tel Aviv 69978, Israel. Email: asafico@tau.ac.il. Supported in part
by ERC Consolidator Grant 863438.} \and Raphael Yuster\thanks{Department of
Mathematics, University of Haifa, Haifa 31905, Israel. E--mail:
raphy@math.haifa.ac.il} }

\date{}

\documentclass [letterpaper,11pt]{article}
\usepackage{amsfonts}

\newtheorem{theo}{Theorem}
\newtheorem{coro}[theo]{Corollary}

\newtheorem{prop}[theo]{Proposition}
\newtheorem{lemma}{Lemma}[section]

\newcommand{\qed}{\hspace*{\fill} \rule{7pt}{7pt}}

\newcommand{\ignore}[1]{}

\begin{document}
\maketitle

\begin{abstract}
The Erd\H{o}s--Stone Theorem asserts that if a graph has edge density $1-1/r+\delta$ then it contains
a complete $(r+1)$-partite graph with $b$ vertices in each part, where $b=b_n(r,\delta) \gg 1$. The celebrated Chv\'atal--Szemer\'edi theorem determined the exact order of $b_n(r,\delta)$ for every $\delta < 1/r^3$. Their bound, however, is not tight when $\delta=1/r-\epsilon$, that is, when the graph has edge density $1-\epsilon$ for small $\epsilon$.
Our main result in this paper determines the correct order in this remaining regime, thereby enabling us to give a tight bound for the Erd\H{o}s--Stone problem for all edge densities. More precisely, we prove that for every integer $r\geq 2$ and $0< \delta < 1/r$ we have
$$
b_n(r,\delta)=\Theta\left(\frac{\log n}{(1/r-\delta)r\log(1/\delta)}\right)\;.
$$
The lower bound is obtained using a K\"{o}vari-S\'os-Tur\'an-type argument combined with a variant of Nikiforov's method of constructing
large blow-ups, while the upper bound is proved using a correlated random graph construction, related to tensor powers.


\end{abstract}

\section{Introduction}\label{intro}

Tur\'an's Theorem \cite{T} asserts (in its weak form) that every graph with more than $(1-1/r)n^2/2$ edges
contains a copy of $K_{r+1}$, the complete graph
on $r+1$ vertices. A $b$-blow-up of $K_{r+1}$, denoted by $K_{r+1}(b)$, is the graph
consisting of $r+1$ disjoint vertex sets $V_1,\ldots,V_{r+1}$ of size $b$ each, where for every $i \neq j$
every vertex in $V_i$ is connected to every vertex in $V_j$.
One of the central results in graph theory, the Erd\H{o}s-Stone Theorem \cite{ErStone}, strengthens Tur\'an's Theorem
by asserting that if $G$ has $(1-1/r+\delta)n^2/2$ edges
then the graph actually contains a blow-up of
$K_{r+1}$. More precisely, there is a function $b_n(r,\delta)$ such that
every $n$-vertex graph with $(1-1/r+\delta)n^2/2$ edges
contains $K_{r+1}(b)$, and such that for fixed $r$ and $\delta$ we have $\lim_{n \rightarrow \infty }b_n(r,\delta)=\infty$.

The determination of the growth rate of $b_n(r,\delta)$ has received considerable attention.
The proof of Erd\H{o}s and Stone \cite{ErStone} showed (roughly) that $b_n(r,\delta)=\Omega(\log^{(r)}n)$,
where $\log^{(r)}x$ denotes the $r$-times iterated logarithm function. This bound was later improved by
Erd\H{o}s \cite{Erdos}. Bollob\'{a}s and Erd\H{o}s \cite{BE} were the first to determine the correct dependence on $n$, showing that
$b_n(r,\delta)=\Theta_{r,\delta}(\log n)$. Bollob\'{a}s, Erd\H{o}s and Simonovits \cite{BES}
further improved the dependence on $\delta$ proving that for large enough $n$, $b_n(r,\delta) \ge c \log n/r \log(1/\delta)$ for some absolute $c > 0$.
They conjectured that $r$ can be removed from the denominator. This conjecture was proved by Chv\'atal and Szemer\'edi \cite{CS}, who showed that for every $r \geq 2$, $0 < \delta \leq 1/r$ and large enough $n$,
\begin{equation}\label{CS}
b_n(r,\delta) \geq  \frac{\log n}{500\log (1/\delta)}\;.
\end{equation}
We finally note that, as proved in \cite{BK,I}, one can also obtain a skewed version in which $r$ of the parts are of logarithmic size,
and the other part is of size polynomial in $n$.

As observed in \cite{BE,BES,CS}, for any fixed $\gamma > 0$, the bound in (\ref{CS}) is tight when $0 < \delta < 1/r^{2+\gamma}$.
In particular, it is tight for, say, $0 < \delta < 1/r^3$.
When $\delta=1/r-\epsilon$, that is, when the graph has $(1-\epsilon)n^2/2$ edges,
we clearly expect to get a bound that increases as $\epsilon \rightarrow 0$. However, note that (\ref{CS}) never gives a bound better than $\log n/\log r$. Nevertheless, it is not hard to use (\ref{CS}) in a black-box fashion\footnote{See the proof of Theorem \ref{theomain} where
we use this observation. We also note that after a close inspection of \cite{CS} we have reached the conclusion that their proof cannot
be modified to give a bound better than $\log n/(r\epsilon\log(1/\epsilon))$. } and show that $b_n(r,1/r-\epsilon) \geq \log n/(r\epsilon\log(1/\epsilon))$. But this turns out not to be the right dependence on $\epsilon$.
The main result of this note provides the following improved bound in this range of parameters.

\begin{theo}\label{theomain}
For every $r \ge 2$, $0 < \epsilon \leq 1/r-1/r^3$ and every $n>n_0(r,\epsilon)$ we have
\begin{equation}\label{theomaineq}
b_n(r,1/r-\epsilon) \ge \frac{\log n}{6000\epsilon r \log r }\;.
\end{equation}
\end{theo}

Needless to say, the constant 6000 in Theorem \ref{theomain} can certainly be improved, at the price of complicating the calculations.
The following proposition shows that the bound in Theorem \ref{theomain} is tight.

\begin{prop}\label{propupper}
For every $r \ge 2$, $0 < \epsilon < 1/r$ and every $n>n_0(r,\epsilon)$ we have
\begin{equation}\label{propeq}
b_n(r,1/r-\epsilon) \leq \frac{10\log n}{\epsilon r \log r }\;.
\end{equation}
\end{prop}

As we mentioned earlier, the Chv\'atal--Szemer\'edi bound (\ref{CS}) is tight in the range
$1/r-1/r^3 < \epsilon < 1/r$. Combining this fact with Theorem \ref{theomain} and Proposition \ref{propupper} shows
that for every $r \ge 2$, $0 < \epsilon < 1/r$ and every $n>n_0(r,\epsilon)$ we have
\begin{equation}\label{coroeq}
b_n(r,1/r-\epsilon) = \Theta\left(\frac{\log n}{\epsilon r \log(r/(1-\epsilon r))}\right)\;.
\end{equation}
We thus obtain a tight bound for the Erd\H{o}s--Stone problem for the entire range of the parameters.
If one prefers the $\delta$ parametrization as in (\ref{CS}), then (\ref{coroeq}) takes the following equivalent form.

\begin{coro}
For every $r \ge 2$, $0 < \delta < 1/r$ and every $n>n_0(r,\delta)$ we have
$$
b_n(r,\delta)=\Theta\left(\frac{\log n}{(1/r-\delta)r\log(1/\delta)}\right)\;.
$$
\end{coro}

The Erd\H{o}s-Stone problem connects the edge density of a graph with the size of the largest blow-up of the complete graph $K_{r+1}$.
A natural related problem asks how the $H$-density of a graph controls the size of its largest $H$-blow-up.
A notable result in this direction is due to Nikiforov \cite{N,N2}, who proved
that an $H$-density of at least $\gamma$ guarantees an $H$-blow-up of size $C_H(\gamma)\log n$. As in the Erd\H{o}s-Stone problem, it is natural
to further optimize the dependence on $\gamma$, a problem studied in \cite{FLW,FWZ,GHW,RS,SY}.

\paragraph{Proof overview:}
Throughout the paper we omit floor/ceiling signs since they are immaterial for our asymptotic estimates.
As opposed to the proof of (\ref{CS}), which relies on an involved application of the regularity lemma,
the proof of Theorem \ref{theomain} is closer in spirit to the proofs of the bounds preceding (\ref{CS}), but with a more modern twist.
It proceeds by constructing the blow-up iteratively, one part at a time. To this end, we rely
on a variant of Nikiforov's inductive approach \cite{N,N2}, used to study the problem we just discussed above. To facilitate this approach, we also require a skewed variant of the K\"{o}vari-S\'os-Tur\'an theorem \cite{KST}.
This proof appears in Section \ref{SecEquiv}.

The construction proving Proposition \ref{propupper} also differs from the classical construction showing the sharpness of (\ref{CS}).
Indeed, the proof of the latter starts with a Tur\'an graph and then adds a random graph inside one of the parts. 
The situation here seems more challenging.
Indeed, one might first take $G(n,1-\epsilon)$ as a candidate example for Proposition \ref{propupper}
but this construction misses (\ref{propeq}) by a $\log r$ factor. To overcome this we use what can be considered a correlated random construction.
We start with a tensor power of the complete graph, thus obtaining a large graph of density $1-\epsilon$ and then sample vertices from it. The main challenge is to pick the parameters so that with high probability the sampled graph contains no large blow-up of $K_{r+1}$. 
This proof appears in Section \ref{SecUpper}.

\paragraph{Declaration of AI usage:} The proof of Proposition \ref{propupper} was suggested to us by ChatGPT. It was not involved in the rest of the paper.

\section{Proof of Theorem \ref{theomain}}\label{SecEquiv}

For the proof of Theorem \ref{theomain} we will need the following K\"{o}vari-S\'os-Tur\'an-type \cite{KST} argument, whose
important point is the following. As is well known, the largest $t$ for which we are guaranteed to find a copy of $K_{t,t}$ in every $n \times n$ bipartite graph with $(1-\zeta)n^2$ edges is $t=\Theta(\frac{\log n}{\zeta})$. The following lemma shows that we can get almost the same bound even if the bipartite graph is skewed to the point of one side being just slightly larger than $t$.

\begin{lemma}\label{KST} The following holds for every $\zeta >0$ and $n \geq n_0(\zeta)$.
Let $F$ be a bipartite graph on vertex sets $A$ and $B$ and
suppose $|A|=n$, $|B|=b$ and $F$ has at least $(1-\zeta)nb$ edges. For a positive integer $d$ satisfying $\zeta \leq 1/(4d)$, set
$t=\frac{\log (n)}{12\zeta\log(d+1)}$. If $b \geq (1+1/d)t$, then $F$ contains a copy of $K_{t,t}$.
\end{lemma}

\paragraph{Proof:} For every vertex $x \in A$ let $d(x)$ denote the degree of $x$. The number of copies
of $K_{1,t}$ with one vertex in $A$ and $t$ vertices in $B$ is
$$
\sum_{x \in A}{d(x) \choose t} \geq n {\frac{1}{n}\sum_{x \in A} d(x) \choose t} \ge n{(1-\zeta)b \choose t}\;,
$$
where the first inequality is Jensen's inequality applied to the convex\footnote{If one wishes to be precise,
then one should use the convex function $f(x)={x \choose t}$ when $x \geq t-1$ and $f(x)=0$ otherwise. Also,
note that the assumptions $b \geq (1+1/d)t$ and $\zeta \leq 1/(4d)$ guarantee that ${(1-\zeta)b \choose t} > 0$.} function ${x \choose t}$,
and the second inequality follows by our assumption on the number of edges of $F$.
Therefore, there must be at least one set $B' \subseteq B$ of $t$ vertices in $B$ which forms a copy of $K_{1,t}$ with at least
$$
n \cdot \frac{{(1-\zeta)b \choose t}}{{b \choose t}} = n \cdot \frac{((1-\zeta)b)!(b-t)!}{((1-\zeta)b-t)!b!}
= n \cdot \prod^{t-1}_{i=0}\left(1-\frac{\zeta b}{b-i}\right)
$$
of the vertices of $A$. To estimate the product on the right hand side, suppose $b=ct$ for some $c > 1$ (recall that we assume that $c \geq 1+1/d$). Then

\begin{eqnarray*}
n \cdot \prod^{t-1}_{i=0}\left(1-\frac{\zeta ct}{ct-i}\right) &\geq& n \cdot \exp\left(-2\zeta ct{\sum^{t-1}_{i=0}\frac{1}{ct-i}}\right) \\
&\geq& n \cdot \exp\left(-2\zeta ct\left(\int^{ct}_{ct-t}\frac{dx}{x}\right)\right) \\
&= &   n \cdot \exp\left(-2\zeta t c\ln\left(\frac{c}{c-1}\right)\right) \\
&\geq & n \cdot \exp\left(-4\zeta t\ln(d+1) \right) \\
&\geq & n^{1/3} \\
&\geq& t \;,
\end{eqnarray*}
where in the first inequality we use the fact that $1-p \geq e^{-2p}$ for every\footnote{The assumption $\zeta \leq 1/(4d)$ guarantees that
$\zeta ct/(ct-i) \leq \zeta c/(c-1) \leq \zeta(d+1) \leq 1/2$.} $0 < p \leq 1/2$, in the third inequality the fact that
$c \geq 1+1/d $ and that the function $f(x)=x\ln(x/(x-1))$ decreases in $(1,\infty)$, and the last inequality relies on the assumption
that $n \geq n_0(\zeta)$. Hence, the set $B'$ along with
$t$ of the vertices of $A$ that form a $K_{1,t}$ with it, form a copy of $K_{t,t}$. $\qed$

We now turn to the main part of this paper, where we prove Theorem \ref{theomain} under a slightly stronger assumption on $\epsilon$.

\begin{lemma}\label{lemmamain}
For every $r \ge 2$, $0 < \epsilon \leq 1/(4r^2)$ and every $n>n_0(r,\epsilon)$ we have
$$
b_n(r,1/r-\epsilon) \ge \frac{\log n}{12\epsilon r \log r }\;.
$$
\end{lemma}

In the proof of Lemma \ref{lemmamain} it will be convenient to assume that the graph has minimum degree $(1-\epsilon)n$ (and not just average
degree $(1-\epsilon)n$). To this end, the following lemma, used (either implicitly or explicitly) in many papers on the Erd\H{o}s-Stone Theorem,
will be useful. For completeness, we give the short proof, similar to the one in \cite{B} p. 330.
\begin{lemma}\label{mindeg}
For every $0 < \eta < 1$ and $n > 4/\eta$ the following holds. If $\epsilon < 1/2$ and $G$ is a graph with $n$ vertices and
$(1-\epsilon){n \choose 2}$ edges, then $G$ has a subgraph $H$ with $p > (\eta/2)^{1/2}n$ vertices and with minimum degree at least
$(1-\epsilon-\eta)p$.
\end{lemma}
\paragraph{Proof:}
Define a sequence of graphs $G_n=G,G_{n-1},G_{n-2},\ldots$ as follows. If the minimum degree of $G_k$ is at least $(1-\epsilon-\eta)k$ then halt.
Otherwise, obtain $G_{k-1}$ from $G_k$ by deleting a minimum degree vertex. Let $H=G_{p}$ be the last graph in the sequence. We have:
\begin{eqnarray*}
\frac{p^2}{2} > {p \choose 2} & \ge & (1-\epsilon){n \choose 2}-(1-\epsilon-\eta)\sum_{k=p+1}^nk\\
              & =   & (1-\epsilon){n \choose 2}-\frac{1-\epsilon-\eta}{2}(n-p)(n+p+1) \\
              &   = & \frac{\eta}{2}(n^2-p^2)-\frac{1-\epsilon-\eta}{2}(n-p)+p^2\frac{1-\epsilon}{2} - \frac{1-\epsilon}{2}n\;.
\end{eqnarray*}
This implies that
$$
p^2(\epsilon+\eta) \ge \eta n^2  - (1-\epsilon-\eta)(n-p)- (1-\epsilon)n \ge \eta n^2-2(1-\epsilon)n\;.
$$
Hence
$$
p^2 \ge \frac{\eta n^2-2(1-\epsilon)n}{\epsilon+\eta} > \frac{\eta}{2}n^2\;,
$$
where in the last inequality we have used $\epsilon < 1/2$ and $n > 4/\eta$.
$\qed$

\paragraph{Proof of Lemma \ref{lemmamain}:}
Let $G$ be an $n$-vertex graph with $(1-\epsilon){n \choose 2}$ edges and suppose
$n$ is large enough so that Lemma \ref{KST} holds whenever $\zeta \geq \epsilon$.
Let $r \ge 2$ be an integer and suppose $\epsilon \le 1/(4r^2)$.
The asymptotic statement of the lemma, together with Lemma \ref{mindeg}, allows us to assume that the minimum degree
of $G$ is at least $(1-\epsilon)n$.

We claim that the graph $G$ contains
a $b_2$-blow-up of $K_2$ ($K_2$ is an edge) with $b_2=\frac{\log n}{12\epsilon \log 2}$ and that for every
$3 \leq i \leq r+1$ the graph $G$ contains
a $b_i$-blow-up of $K_i$ with $b_i=\frac{\log n}{12\epsilon (i-1)\log (i-1)}$. We proceed by induction on $i$. For the case $i=2$
let us define a bipartite graph on vertex sets $A$ and $B$, where each of these sets contains $n$ vertices and vertex
$a_i \in A$ is connected to $b_j \in B$ if and only if $(i,j)$ is an edge of $G$. Setting $t=b_2$,
$d=1$, and $\zeta = \epsilon$ we see that since $n \gg 2t$, we can conclude using Lemma \ref{KST} that $F$ contains a copy of $K_{b_2,b_2}$ which corresponds to a $K_{b_2,b_2}$ in $G$. Assuming we have established the result for $i-1 \leq r$, we prove it for $i$.

Let $K$ be a $b_{i-1}$-blow-up of $K_{i-1}$ whose existence is guaranteed
by the induction hypothesis. We now show how to find in $G$ a $b_i$-blow-up of $K_{i}$.
Partition the vertices of $K$ into
$b_{i-1}$ vertex-disjoint copies of $K_{i-1}$ and denote them by $S_1,\ldots,S_{b_{i-1}}$.
Let $F$ be the bipartite graph on vertex sets $A$ and $B$,
where the vertices of $A$ are the $n$ vertices of $G$, and $B$ has $b_{i-1}$ vertices representing $S_1,\ldots,S_{b_{i-1}}$. We put an
edge in $F$ between $a \in A$ and $S_j \in B$ if and only if vertex $a$ is connected in $G$ to all the vertices of $S_j$. Since the
minimum degree of $G$ is $(1-\epsilon)n$, we get that every vertex in $B$ is connected to at least $n(1-\epsilon(i-1))$ of the
vertices of $A$. In particular, $F$ has at least $(1-\epsilon(i-1))nb_{i-1}$ edges.
Now observe that for $i\geq 4$
$$
b_{i-1}=\frac{\log n}{12\epsilon(i-2)\log(i-2)} \ge \left(1+1/(i-2)\right) \frac{\log n}{12\epsilon (i-1)\log(i-1)} = \left(1+1/(i-2)\right)b_{i}\;.
$$
The same also holds for $i=3$ since $b_2=2b_3$.
Hence, for every $i \geq 3$ we can apply Lemma \ref{KST} on $F$, with $\zeta=\epsilon(i-1)$, $d=i-2$, $b=b_{i-1}$ and $t=b_i$, to conclude that it contains a copy of $K_{b_i,b_i}$. Note that the lemma's assumption that $\epsilon \le 1/(4r^2)$ guarantees that $\zeta \leq 1/(4d)$ as needed for the application of Lemma \ref{KST}.
Suppose the vertices of this copy are $a_1,\ldots,a_{b_i} \in A$ and $S_1,\ldots,S_{b_i} \in B$. It is now easy to see that the definition of $F$ guarantees that the vertices
$a_1,\ldots,a_{b_i}$ along with the $(i-1)b_{i}$ vertices of $\bigcup^{b_i}_{j=1}S_j$ form a $b_i$-blow-up of $K_i$ in $G$. This completes the induction step and the lemma follows from the case $i=r+1$. $\qed$

\bigskip

To prove Theorem \ref{theomain}, it remains to treat the range of $\epsilon$ not covered by Lemma \ref{lemmamain}.
We will do so by reduction to (\ref{CS}).

\paragraph{Proof of Theorem \ref{theomain}:} If $\epsilon \leq 1/(4r^2)$, then the result follows from Lemma \ref{lemmamain}. Suppose now that
$1/(4r) \leq \epsilon \leq 1/r -1/r^3$, then we get from (\ref{CS}), with $\delta=1/r^3$, that $G$ has $b$-blow-up of $K_{r+1}$ with
$b=\log n/(1500\log r)$ which agrees with (\ref{theomaineq}) when $1/(4r) \leq \epsilon \leq 1/r -1/r^3$.

So suppose for the rest of the proof that $1/(4r^2) \leq \epsilon \leq 1/(4r)$. Applying (\ref{CS}) with $r'=\lfloor 1/2\epsilon \rfloor$ and $\delta=\epsilon$ we get that
$G$ contains a $b'$-blow-up of $K_{r'+1}$, with
$b' = \log n/500\log (1/\epsilon)$.
We can now group the $r'+1$ partition classes of
this blow-up into $r+1$ sets\footnote{This is indeed possible since we assume that $\epsilon \leq 1/4r$.}
to get a $b$-blow-up of $K_{r+1}$ with
$$
b= \frac{b'(r'+1)}{(r+1)}  \geq \frac{\log n}{(8/3)500\epsilon r \log 1/\epsilon} > \frac{\log n}{6000\epsilon r \log r}\;,
$$
where the last inequality uses the assumption $\epsilon \geq 1/(4r^2)$. $\qed$

\bigskip

\section{Proof of Proposition \ref{propupper}}\label{SecUpper}

Set $b=\frac{10\log n}{\epsilon r \log r}$, $m=(r+1)b$, $q=\frac{m}{2}$ and $\ell=\frac{\epsilon q}{2}=\frac{\epsilon m}{4}$.
Let $H=K^{\otimes \ell}_q$ denote the $\ell$-th tensor power of the complete graph $K_q$, namely, the vertices of $H$
are $[q]^{\ell}$ and two vertices $x=(x_1,\ldots,x_{\ell})$ and $y=(y_1,\ldots,y_{\ell})$ are connected if and only if $x_i \neq y_i$ for every $1 \leq i \leq \ell$. To construct the $n$-vertex graph $G$ on vertex set $[n]$ establishing Proposition \ref{propupper} we will sample $n$ vertices $v_1,\ldots,v_n$ from $H$ with repetitions, and put $(i,j) \in E(G)$ if and only if $(v_i,v_j) \in E(H)$. We will now show that with probability at least $3/4$ the edge density of $G$ is at least $1-\epsilon$, and that with probability at least $3/4$ the graph $G$ has
no $b$-blow-up of $K_{r+1}$. These two assertions would imply the upper bound on $b_n(r,1/r-\epsilon)$.
We will frequently use the fact that each vertex in $[n]$ is assigned a random vertex in $H$, hence we can think of each vertex in $G$ as being assigned a random string in $[q]^{\ell}$, with all coordinates over all vertex/string choices being independent.

The first assertion is simple. The probability that $(i,j) \in E(G)$ is $(1-1/q)^{\ell} \geq 1-\ell/q \geq 1-\epsilon/2$.
We can imagine the pairs of vertices of $G$ as $n$ matchings of size $n/2$ (assuming, as we may, that $n$ is even).
Then in each matching we can apply a standard additive Chernoff bound to get that the fraction of non-edges is at most $\epsilon$
with probability $1-e^{-\frac14\epsilon^2n}$. Since $n > n_0(\epsilon)$, a union bound over all $n$ matchings, gives that $G$ has
edge density at least $1-\epsilon$ with probability at least $1-ne^{-\frac14\epsilon^2n} \geq 3/4$.

We now prove the second assertion. We will prove that for every fixed collection of $r+1$ disjoint vertex sets $B_1,\ldots,B_{r+1}$ in $[n]$ of size $b$ each, the probability that they form a copy of the $b$-blow-up of $K_{r+1}$ is at most $(r+1)^{-\epsilon m^2/8}$. Taking a union bound over all $n^m$ possible choices of $B_1,\ldots,B_{r+1}$ and recalling that $m=\frac{10(r+1)\log n}{\epsilon r \log r}$ would give that $G$ satisfies the second assertion with probability at least $1-n^m(r+1)^{-\epsilon m^2/8} \geq 3/4$.

So let us then fix disjoint $B_1,\ldots,B_{r+1}$ as above. For every $1 \leq i \leq r+1$ and $1 \leq j \leq \ell$
let $D_{i,j} \subseteq [q]$ denote the numbers in $[q]$ used in the $j$-{th} coordinate of at least one vertex/string of $B_i$. The crucial
observation now is that if $B_1,\ldots,B_{r+1}$ form a blow-up of $K_{r+1}$ then for every $1 \leq j \leq \ell$, the sets $D_{1,j},\ldots,D_{r+1,j}$ are disjoint. Let us estimate the probability of this event for some fixed $1 \leq j \leq \ell$.

There are $(r+1)^q$ ways to pick disjoint sets $D_1,\ldots,D_{r+1}$ whose union is $[q]$.
Fix one such choice. The probability that $D_{i,j} \subseteq D_{i}$ is $(|D_i|/q)^b$, so the probability
this happens for all $1 \leq i \leq r+1$ is $(|D_1|/q)^b \cdots (D_{r+1}/q)^b$. Since $D_{1},\ldots,D_{r+1}$ are disjoint
we have $|D_1|+\ldots+|D_{r+1}| = q$. Hence, by Jensen's inequality we have
$$
\left(\frac{|D_1|}{q}\right)^b \cdots \left(\frac{|D_{r+1|}}{q}\right)^b \leq \left(\frac{1}{r+1}\right)^{(r+1)b}=\left(\frac{1}{r+1}\right)^{m}\;.
$$
We conclude that the probability that the sets $D_{1,j},\ldots,D_{r+1,j}$ are disjoint is at most $(r+1)^{q-m}=(r+1)^{-m/2}$

Since the entries of the vertices/strings are independent, the probability that
$D_{1,j},\ldots,D_{r+1,j}$ are disjoint is independent of the probability that $D_{1,j'},\ldots,D_{r+1,j'}$ are disjoint,
hence the probability that $B_1,\ldots,B_{r+1}$ form a blow-up of $K_{r+1}$ is at most $[(r+1)^{-m/2}]^{\ell}=(r+1)^{-\epsilon m^2/8}$,
and the proof is complete.


\begin{thebibliography}{99}

\bibitem{B}
B. Bollob\'{a}s, {\bf Extremal Graph Theory}, Academic Press, London 1978.

\bibitem{BE}
B. Bollob\'{a}s and P. Erd\H{o}s,
{\em On the structure of edge graphs},
Bull. London Math. Soc. 5 (1973), 317--321.

\bibitem{BES}
B. Bollob\'{a}s, P. Erd\H{o}s, and M. Simonovits,
{\em On the structure of edge graphs II},
J. London Math. Soc. 12 (1976), 219--224.

\bibitem{BK}
B. Bollob\'{a}s and Y. Kohayakawa,
{\em An extension of the Erd\H{o}s-Stone Theorem},
Combinatorica 14 (1994), 279--286.

\bibitem{CS}
V. Chv\'atal and E. Szemer\'edi,
{\em On the Erd\H{o}s-Stone Theorem}.
J. London Math. Soc. 23 (1981), 207--214.

\bibitem{Erdos}
P. Erd\H{o}s, Some recent results on extremal problems in graph theory, in ``Actes des
journees detudes sur la theorie des graphes,'' pp. 117--130, Dunod, Paris, 1967.


\bibitem{ErStone}
P. Erd\H{o}s and A.~H. Stone,
{\em On the structure of linear graphs},
Bull. Amer. Math. Soc. 52 (1946), 1089--1091.

\bibitem{FLW}
J. Fox, S. Luo and Y. Wigderson, Extremal and Ramsey results on graph blowups, J.
Comb. 12 (2021), 1--15.

\bibitem{FWZ}
J. Fox, Y. Wigderson and Y. Zhou, Finding blowups one vertex at a time, arXiv:2605.23301, 2026.


\bibitem{GHW}
A. Girao, Z. Hunter and Y. Wigderson, Blowups of triangle-free graphs, Advances in Combinatorics 10, 2025.


\bibitem{I}
Y. Ishigami, Proof of a conjecture of Bollob\'{a}s and Kohayakawa on the Erd\H{o}s-Stone Theorem,
J. Combin. Theory Ser. B 85 (2002), 222--254.

\bibitem{KST}
T. K\"{o}vari, V. S\'os, and P. Tur\'an,
{\em On a problem of K. Zarankiewicz},
Colloquium Math. 3 (1954), 50--57.

\bibitem{N}
V. Nikiforov, {\em Graphs with many $r$-cliques have large complete $r$-partite subgraphs},
Bull. London Math. Soc. 40 (2008), 23--25.

\bibitem{N2}
V. Nikiforov, Graphs with many copies of a given subgraph, Electronic Journal of Combinatorics 15 (2008), Note N6.

\bibitem{RS}
V. R\"odl and M. Schacht, Complete partite subgraphs in dense hypergraphs, Random
Structures Algorithms 41 (2012), 557--573.

\bibitem{SY}
A. Shapira and R. Yuster, On the Density of a Graph and its Blowup, J. Combin. Theory Ser. B 100 (2010), 704--719.

\bibitem{T}
P. Tur\'an,
{\em On an extremal problem in graph theory} (In Hungarian),
Mat. Fiz. Lapok 48 (1941), 436--452.

\end{thebibliography}
\end{document}